\documentclass[12pt]{article}
\usepackage{amssymb}

\newcommand\R{{\mathbb R}}
\newcommand\C{{\mathbb C}} 
\newcommand\Z{{\mathbb Z}}
\newcommand\Q{{\mathbb Q}}
\newcommand\A{{\mathbb A}}
\newcommand\B{{\cal B}}
\renewcommand\a{{\cal A}}
\newcommand\G{{\Gamma}}
\newcommand\g{{\gamma}} 
\newcommand\gl{{\bf g}} 
\newcommand\tg{{\widetilde{\g}}}
\newcommand\gm{{\bf G}_m} 
\renewcommand\d{{\partial}} 
\renewcommand\P{{\mathbb P}} 
\renewcommand\k{{\bf k}} 
\newcommand\ra{{\rightarrow}}
\newcommand\1{{\bf 1}}
\newcommand\mod{{\mbox {mod}}}

\newcommand\Cech{{\v{C}ech }}
\newcommand\E{{\cal E}}
\newcommand\h{{\hbar}} 

\newcommand\jt{{jets_\infty(T_{poly})}}
\newcommand\jd{{jets_\infty(D_{poly})}}
\newcommand\qed{\,\,{$\blacksquare$}}

\renewcommand\O{{\cal O}}

\renewcommand\L{{\cal L}}

\newtheorem{lemma}{Lemma}
\newtheorem{rmk}{Remark}
\newtheorem{prop}{Proposition}
\newtheorem{dfn}{Definition} 
 \newtheorem{conj}{Conjecture} 
\newtheorem{const}{Construction}
\newtheorem{ques}{Question}
\newtheorem{quess}[ques]{Questions}

\begin{document}

\title{Deformation quantization of algebraic varieties}

\author{Maxim Kontsevich\\
{\small \textit {IH\'ES, 35, route de Chartres, F-91440 Bures-sur-Yvette, 
France}}\\  {\small email: {\texttt{maxim@ihes.fr}}}}

\maketitle

\begin{abstract}
The paper is devoted to peculiarities of the deformation quantization
in the algebro-geometric context. A direct application of the formality 
theorem to an algebraic Poisson manifold gives a
canonical sheaf of categories deforming coherent sheaves.
The global category is very degenerate in general.
Thus, we introduce a new notion of a semi-formal deformation, a replacement 
in algebraic geometry of an actual deformation (versus a formal one).
Deformed algebras obtained by semi-formal deformations are Noetherian
and have polynomial growth.

We propose constructions of semi-formal quantizations of
projective and affine algebraic Poisson
manifolds satisfying certain natural geometric conditions.
Projective symplectic manifolds (e.g. K3 surfaces and abelian varieties)
do not satisfy our conditions, but projective spaces with quadratic 
Poisson brackets and Poisson-Lie groups can be semi-formally quantized.
\end {abstract}

\section*{Introduction} 

In the present paper we describe applications in the algebro-geometric
setting of the formality theorem from \cite{k1} (where it is proven in the
$C^\infty$ situation).  The general question is whether it is possible to
``quantize'' (i.e. to deform forgetting commutativity) an algebraic Poisson
manifold, and obtain a kind of a noncommutative space.  

One can consider schemes in ``noncommutative algebraic geometry in
quasiclassical regime'' as spaces endowed with sheaves of 
 algebras over the ring of formal power series.  In fact, the
formal deformation quantization gives a little  bit more complicated
structure, a sheaf (or, in another terminology, a stack) of algebroids. In
Section 1 we describe briefly this structure, and in Appendix A we give 
directions towards explicit formulas. These formulas are quite complicated;
in a sense we use the formality morphism from \cite{k1} as a ``black box''. 
  
Using the stack of algebroids we define a canonical abelian category of
coherent sheaves on the quantized manifold. Technical tools involved in the
definition  are borrowed from the theory of usual (commutative) formal
noetherian schemes. In a sense, a formal quantization of an algebraic Poisson
manifold can be seen as  a noncommutative version of the usual
formal thickening. 

 For applications it is more interesting to have  not only formal  but actual
deformations. At the end of this Introduction we  introduce a  notion
of a {\em semi-formal} deformation of an algebra, which is a reasonably good 
substitute in algebraic geometry for the notion of an actual  deformation
(versus a purely formal one).

In Section 2 we discuss examples of quantizations of the affine space.
 It turns out that  in general one can quantize semi-formally only
 Poisson brackets whose coefficients are at most quadratic polynomials. The
canonical correspondence between classical and quantum quadratic algebras give
rise to an interesting (formal) transcendental map between two moduli schemes.

In Sections 3 and 4 we give simple geometric criteria for projective 
(resp. affine) Poisson manifolds to admit a semi-formal noncommutative 
deformation over the ring of formal power series. Our criteria of 
quantizability seem to be not too far from  optimal ones.
After localization to the field of Laurent series (or assuming convergence
 of the perturbation series in p-adic or archimedean topology) we obtain
objects which can be considered as genuine noncommutative schemes.
It turns out that abelian varieties with non-zero Poisson structures
admit {\em no} semi-formal quantization, in contrast with the $C^\infty$
situation where quantum tori are important examples.
  
In Appendix we describe a general construction of algebroids, and give a
 parametrization of formality morphisms in the global $C^\infty$ case by
connections on the tangent bundle of the underlying manifold.

\subsection*{A convention concerning associators}
Explicit formulas for the formality morphism have real (resp. complex)
coefficients, the weights of graphs in~\cite{k1}, (resp. coefficients in 
the Drinfeld associator in~\cite{t}). These coefficients satisfy an 
infinite system of non-homogeneous  quadratic equations defined 
over $\Q$.  One can show (see e.g ~\cite{k2} and \cite{d1}) that in
both approaches \cite{k1} and \cite{t} there exists a rational solution.
In fact, different solutions can give  \emph{non-isomorphic} 
deformation quantizations procedures, at least when one not only
quantizes plain Poisson manifolds but considers the formality morphism 
in all degrees (see Sect. 4.6. in \cite{k2}).  

In the present paper, when we consider algebraic varieties over a general 
field $\k$ of characteristic zero, we always assume that a solution of the 
system of quadratic equations as above is given. Strictly speaking, 
true parameters for the  formal deformation quantization are  pairs 
(a Poisson manifold, a choice of the Drinfeld associator).

\subsection*{Formal deformations, compactness and filtrations}

Formal deformation theory applied to ``non-compact'' objects usually 
gives nonsensical results. For example, if $(X_t)_{t\in C}$ is a smooth 
one-parameter  family of smooth affine varieties over an algebraic curve $C$, 
then  after passing to completion at any point $t_0\in C$  we always 
get a {\em formally trivial} family.   

In the case of associative algebras the role of compact objects is 
played by finite-dimensional algebras.  Let $R$ be  a complete pro-Artin 
local ring with residue field $\k$. 
By definition, a deformation over $Spec(R)$ of a $\k$-algebra $A$ 
is an algebra $\widehat{A}$ over $R$, {\em topologically free} 
as $R$-module, considered together  with an identification of 
$\k$-algebras $A\simeq\widehat{A}\otimes_R \k$.    
Notice that the algebra $\widehat{A}$ is free as $R$-module if and only if 
$A$ is finite-dimensional and $R$  has finite length.  

There is still a way to ``compactify'' infinite-dimensional algebras.
 Namely, if an algebra $A=\oplus_{n\in \Z} A^n$ is $\Z$-graded and finitely
generated, and  all graded components  $A^n$ are finite-dimensional, 
then deformations of $A$ considered as a {\em graded} algebra are similar
to deformations of a finite-dimensional algebra.
The deformed graded algebra is always free as $R$-module,
and deformations depend on a finite number of parameters:
there exists a versal deformation (i.e. any other deformation is induced 
from it) over a complete local ring with finitely generated maximal ideal.

The case of a non-graded algebra can be treated similarly.
We assume that the algebra $A$ is finitely generated. 
The following lemma is obvious: 
\begin{lemma} An algebra $A$ over $\k$ is finitely generated iff
 there  exists an increasing filtration of $A$ by finite-dimensional 
subspaces $(A_{\le n})_{n\ge 0}$  such that 
$A_{\le n}\cdot A_{\le m}\subset A_{\le n+m}$, $A=\cup_{n\ge 0} A_{\le n}$ 
and  the associated graded algebra 
$$gr(A):=\oplus_{n\ge 0} A_{\le n}/A_{\le n-1},\,\,A_{\le -1}:=0$$
is finitely generated.
\end{lemma}
  
Let us choose  a filtration of $A$ as in the lemma. 
It is easy to see that the Rees ring   
$$A':=\oplus_{n\ge 0} A_{\le n} t^n\subset A[t]$$ 
is also a finitely generated graded algebra over $\k$. 
Let us assume that we have a deformation $\widehat{A'}$ over $R$ of a 
{\em graded} algebra $A'$, together with a central element $t'$ of
degree $+1$ reducing to $t$ after specialization to the closed point 
$Spec(\k)$ of $Spec(R)$.  Let us denote by $\widehat{A}_{finite}$ the 
degree $0$ component of the localized algebra $A'[(t')^{-1}]$. 
It is easy to see that $\widehat{A}_{finite}$ is a dense subalgebra in 
a formal deformation $\widehat{A}$ of $A$, and the $R$-module $\widehat{A}$ 
coincides with the completion  of $\widehat{A}_{finite}$ in the adic topology.
Moreover, $\widehat{A}_{finite}$ is endowed with an increasing  filtration by
finitely generated free $R$-modules, reducing to the filtration  $(A_{\le
n})_{n\ge 0}$ after the specialization. We call such an algebra 
$\widehat{A}_{finite}$ a semi-formal deformation of $A$.  
Here is the precise definition:

\begin{dfn} A {\em semi-formal deformation} over $R$ of a finitely generated
associative algebra $A/\k$ is  an algebra $\widehat{A}_{finite}$ over $R$  
endowed with an identification of $\k$-algebras 
$A\simeq \widehat{A}_{finite}\otimes_R \k$ such that
there exists an exhaustive increasing filtration 
 of $\widehat{A}_{finite}$  by finitely generated
free $R$-modules $\widehat{A}_{\le n}$, compatible
with the product,  admitting a splitting as a filtration of $R$-modules, 
and such that the Rees ring of the induced filtration of $A$ is 
finitely generated over $\k$. 
\end{dfn}

In Section 4 under certain assumptions, we shall construct semi-formal 
deformations of algebras of functions on affine Poisson  manifolds.  

\subsection*{Acknowledgement} I am grateful to A. Rosenberg and to  Y.
Soibelman for very useful discussions.

\section{Formal quantization and algebroids}

\subsection{Algebroids} 

Any connected non-empty groupoid (i.e. a non-empty category 
whose morphisms are invertible and all objects are isomorphic to each other)
is in a   sense a ``group defined up to an inner automorphism''.  Here we
introduce   an analogous notion for associative algebras.  
\begin{dfn} Let $R$ be a commutative ring. 
An \emph{algebroid} over $R$ is a small  category $\a$ such that 
\begin{itemize} 
\item
$R$ is nonempty and all objects of $R$ are isomorphic,
\item morphism sets are endowed with the
structure of $R$-modules,  \item compositions of morphisms are $R$-bilinear,
\end{itemize}
\end{dfn} 

An algebroid over $R$ with one object is the same as a unital associative 
algebra over $R$. Conversely, we can associate with any 
algebroid $\a$ a canonical {\em isomorphism class} of algebras. 
Namely, all algebras $\a_x:=Hom_\a(x,x)$ where $x\in Ob(\a)$  
are isomorphic to each other. For any  algebroid $\a$ we define abelian 
$R$-linear category   $\a-\mod$ as the category of $R$-linear 
functors from $\a$ to the category $R-\mod$ (of left $R$-modules). 
 This category is canonically equivalent to the category 
 $\a_x-\mod$ for any $x\in Ob(\a)$.

 \subsection{Formal quantization of affine Poisson manifolds} 
 Let $X/\k$ be a smooth affine algebraic variety where $\k,\,\,char(\k)=0$ is
a field, and $\g\in \Gamma(X,\wedge^2 T_X),\,\,[\g,\g]=0$ be an algebraic
Poisson  structure on $X$. Our goal is to quantize the Poisson algebra   
$\O(X)$ in the formal sense, i.e. to construct a noncommutative associative   
$\k[[\hbar]]$-linear product on $\O(X)[[\hbar]]$.

In order to quantize  $\O(X)$ one needs to make an additional choice in order
to rigidify the setting. For example, a choice of a connection $\nabla$ on
the tangent bundle to $X$ is sufficient. Such a connection exists always, 
just because the obstruction to its existence lies in 
$H^1(X,T^*_X\otimes_{\O_X} End(T_X))$, 
and the latter space vanishes for affine $X$.

 \begin{const} 
 We  associate canonically with $(X,\gamma)$ an  algebroid  $\a_{X,\gamma}$
over $\k[[\hbar]]$.  Objects of $\a_{X,\gamma}$ are connections on $T_X$.     
Morphisms spaces $Hom_{\a_{X,\gamma}}(\nabla_1,\nabla_2)$ (for any two 
objects $\nabla_1,\nabla_2$) are canonically identified as 
$\k[[\hbar]]$-modules with  $\O(X)[[\hbar]]$.      
Compositions of morphisms 
$$Hom_{\a_{X,\gamma}}(\nabla_1,\nabla_2)\otimes 
Hom_{\a_{X,\gamma}}(\nabla_2,\nabla_3)\ra 
Hom_{\a_{X,\gamma}}(\nabla_1,\nabla_3)$$ 
are given (in local coordinates) by an explicit series in $\hbar$ whose terms 
are bidifferential operators on $X$ with coefficients which are differential 
polynomials in $\nabla_1,\nabla_2,\nabla_3$ and $\gamma$. 
These operators do not depend on the choice of local coordinates. 
 \end{const}

 We shall not present here the above mentioned explicit formulas but indicate
in the Appendix a way to obtain them. 
 
Thus, similarly to the $C^\infty$ case, we associate with $(X,\gamma)$ 
a canonical isomorphism class of algebras. Given a connection $\nabla$ on 
$T_X$ we denote by $\O_{\gamma,\nabla}(X)$ the algebra over
$\k[[\hbar]]$ corresponding to the object $\nabla\in Ob(\a_{X,\gamma})$.
In the particular case of a connection with  both the curvature and the 
torsion vanishing, the product on  $\O_{\gamma,\nabla}(X)$ is given by 
formulas in local formal affine coordinates identical to ones in [K1].
  
   Also, we associate canonically with $(X,\gamma)$ an abelian
$\k[[\hbar]]$-linear category. Recall that for any affine scheme $S$ the
category $QCoh(S)$ of quasi-coherent sheaves on $S$ is canonically 
equivalent to the category $\O(S)-\mod$. 
We denote the canonical abelian category of $\k[[\hbar]]$-linear 
functors from $\a_{X,\gamma}$ to $\k[[\hbar]]-\mod$ simply by 
$QCoh_\gamma(X)$, and consider it as a ``formal quantization'' of  
$QCoh(X)$. Thus, even without the choice of a connection we still can 
speak about ``modules over the formally quantized algebra''.  
Let $f:(X_1,\gamma_1)\ra (X_2,\gamma_2)$ be a Poisson map which is \'etale, 
i.e. the tangent map $Tf_x:T_x X_1\ra T_{f(x)} X_2$ is an isomorphism 
for any $x\in X_1(\overline{\k})$ (an analog of a local diffeomorphism 
in the $C^\infty$ case). 
It follows from the local nature of the formulas for the composition 
of morphisms, that $f$ induces a functor between algebroids
$f^*_{alg}:\a_{X_2,\gamma_2}\ra \a_{X_1,\gamma_1}$. 
Composition with $f^*_{alg}$ gives a direct image functor
$f_*:QCoh_{\gamma_1}(X_1)\ra QCoh_{\gamma_2}(X_2)$.  The left adjoint 
to it is denoted by $f^*$  and called the pullback functor. In particular,
we get restriction functors for open embeddings.  

\begin{rmk} It is easy to see that the construction of the algebroid
 $\a_{X,\gamma}$ extends to arbitrary, not necessarily affine Poisson 
manifolds $(X,\gamma)$. Objects of $\a_{X,\gamma}$ are by definition 
connections on $T_X$.  For general $X$, the algebroid $\a_{X,\gamma}$ 
is always connected (but possibly empty).
\end{rmk}

\subsection{Noetherian property}

It is easy to see that the algebra $\O_{\gamma,\nabla}(X)$
is (left and right) noetherian. It is a corollary of the following  
well-known result (\cite{bou}, Corollary 2 to Proposition 12, 
Chapter III, \S 2, no. 9):

\begin{prop} \label{Noether} A complete  Hausdorff filtered ring with  
an exhaustive filtration is left noetherian if the associated graded ring 
is left noetherian.
\end{prop} 

We denote by $Coh_\gamma(X)$ the full {\em abelian} subcategory of
$QCoh_\gamma(X)$  consisting of functors for which the corresponding modules
over algebras  $\O_{\gamma,\nabla}(X)$ are finitely generated. 

The next Proposition gives a possibility to describe category $Coh_\g(X)$ in
terms of projective limits. 

\begin{prop} \label{fin-inf} Let $A$ be a ring with central element $\h\in A$
such that $A$ coincides with its completion with respect to the two-sided 
ideals generated by powers of $\h$:
$$A=\lim_{\longleftarrow} A/\h^n A$$
and the associated graded ring $gr(A):=\oplus_{n\ge 0} \h^n A/\h^{n+1} A$ 
is noetherian.

Then the  category of finitely  generated $A$-modules is
naturally equivalent to the category  of collections $(M_n,\phi_{n})_{n\ge 1}$
 where $M_n$ is a  module over $A_n:=A/\h^n A$ and 
$\phi_{n}:A_n\otimes_{A_{n+1}} M_{n+1}\ra M_n$ is an isomorphism of 
$A_{n}$-modules, and $A_1$-module $M_1$ is finitely generated.
The correspondence between finitely  generated $A$-modules $M$ and collections 
$(M_n,\phi_{n})_{n\ge 1}$ is the following: 
$$M \mapsto  (M_n:=M/\h^n M=A_n\otimes_A M, \mbox{an obvious epimorphism}
\phi_n)_{n\ge 1}$$ 
$$(M_n,\phi_{n})_{n\ge 1}\mapsto
M:=\lim_{\stackrel{\longleftarrow}{(\phi_n)}} M_n$$ 
\end{prop}

{\it Proof:} We shall check the equivalence only at the level of objects;
for the morphisms the proof of the equivalence is similar.

Let us start with a finitely generated $A$-module $M$. Then
collection $(M/\h^n M)_{n\ge 1}$ with obvious epimorphisms satisfies the
conditions formulated in the Proposition.  We claim that we have an equality 
$$M=\widehat{M}:=\lim_{\longleftarrow} M/\h^n M$$  Let us represent $M$ as a
quotient $M=A^I/B$  where $B\subset A^I$ is a left submodule
in a free finitely generated module $A^I$. The property $M=\widehat{M}$ 
is equivalent to the completeness of $B$ with respect to the filtration 
$B_{\ge n}:=B\cap \h^n A^I$. By the noetherian property of $gr(A)$, the module
$gr(B)\subset gr(A)^I$ is finitely generated. Let us choose a finite system
of homogeneous generators $(x_j\in gr_{d_j}(B))_{j\in J}$ of $gr(B)$ and their
lifts $(\widetilde{x}_j\in B_{\ge d_j})_{j\in J}$. We endow the free module
$A^J$ with the complete $\Z$-filtration as
$$(A^J)_{\ge n}:=\{(a_j)_{j\in J}|\,\forall j \, a_j\in \h^{max(0,n-d_j)}A\}$$
The map $A^J\ra B$ sending canonical generators to elements
$(\widetilde{x}_j)_{j\in J}$ is compatible with filtrations, and the associated
graded map $gr(A^J)\ra gr(B)$ is an epimorphism. By completeness of $A$ we
conclude that  $B$ is also complete.

Now we start with a collection $(M_n,\phi_{n})_{n\ge 1}$ satisfying
conditions as in the Proposition. We claim that 
$M:=\lim_{\longleftarrow} M_n$  is a finitely generated $A$-module 
and the map $M/\h^n M\ra M_n$ is an isomorphism for any $n$.
Let us endow $M$ with the filtration $M_{\ge n}:=Ker(M\ra M_n)$.  
The associated graded module 
$$gr(M)=\oplus_{n\ge 1} Ker(M_{n+1}\ra M_n)$$
is generated by $M_1$, thus it is finitely generated over $gr(A)$.
Therefore $M$ is also finitely generated.
The equality $M_{\ge n}=\h^n M$ follows from the definition of $M$.
\qed

\subsection{Gluing of coherent sheaves}

Let $(X,\g)$ be a smooth affine Poisson algebraic variety over $\k$
and $\nabla$ be a connection on $T_X$. 

Here we shall mention several results which later will allow us to define a
stack of categories  of coherent sheaves. The proofs are very easy and 
left to the reader as an exercise.

 \begin{prop} For any open affine subscheme $U\subset X$, the 
algebra $\O_{\g,\nabla}(U)/(\h^n)$ is flat as right
$\O_{\g,\nabla}(X)/(\h^n)$-module. \end{prop}

Thus, the pullback functor 
$\O_{\g,\nabla}(U)/(\h)^n \otimes_{\O_{\g,\nabla}(X)/(\h)^n}\cdot$
is exact. Obviously it preserves the class of finitely generated modules.
Passing to the limit as $n\ra \infty$ and using Proposition~\ref{fin-inf}
we obtain the exactness of the  pullback functor 
$Coh_\gamma(X)\ra Coh_\gamma(U)$.

\begin{prop} For any two open affine subschemes $U_1,U_2\subset X$
 we have a natural equivalence of bimodules over $\O_{\g,\nabla}(X)/(\h^n)$:
$$\O_{\g,\nabla}(U_1\cap U_2)/(\h^n)=\O_{\g,\nabla}(U_1)/(\h^n)
\otimes_{\O_{\g,\nabla}(X)/(\h^n)}\O_{\g,\nabla}(U_2)/(\h^n)$$
\end{prop}

\begin{prop} If $(U_i)_{i\in I}$ is a covering of $X$ by 
open affine subschemes, then $\oplus_i \O_{\g,\nabla}(U_i)/(\h^n)$ is
 strictly flat as right $\O_{\g,\nabla}(X)/(\h^n)$-module.
\end{prop}
Recall that a flat right module is strictly flat iff the tensor product
with it is a faithful functor.

\subsection{Formal quantization of general Poisson manifolds} 
 
Let $(X,\gamma)$ be  a not  necessarily affine algebraic Poisson manifold.  
It follows from the previous discussion that for any $n\ge 1$ we get 
 a canonical presheaf of algebroids $U\mapsto\a_{U,\gamma}/(\h^n)$
on $X$ in Zariski (or \'etale) topology. The associated sheaf of 
categories (or a stack in the terminology of Giraud)  is denoted  by 
$\a'_{\gamma,n}$. For open affine $U\subset X$, the category
 $\a'_{\gamma,n}(U)$ contains a subalgebroid canonically equivalent to 
$\a_{U,\gamma}/(\h^n)$. Passing to the  projective limit $n\ra \infty$ 
we obtain a stack $\a'_{\gamma}$, which we consider
 as the formal deformation quantization of $(X,\gamma)$.

 It is reasonable to call $\a'_{\gamma}$ a {\em ``stack of algebroids''},
because for any $x\in X$ the stalk ${\left(\a'_{\gamma}\right)}_x$  is an
algebroid.  Categories $\a'_{\gamma}(U)$ for open $U\subset X$ are not
algebroids in general. For example, for the trivial Poisson structure $\g=0$,
the category  $\a'_{\gamma,1}(U)$ is naturally equivalent to the full 
subcategory of $Coh(U)$ consisting of all line bundles on $U$.
 Also, for non-trivial $\gamma$, the global category $\a'_{\gamma}(X)$  
can be empty or can have more than one isomorphism class of objects. 
 
 Similarly, as follows from the propositions in the previous  subsection, 
we can  define a canonical category of coherent sheaves 
$Coh_\gamma(X)/(\h^n)$ using gluing data\footnote{Here we have a particular
example of the gluing in noncommutative algebraic geometry; 
see \cite{kr} for the general approach.}. 
Passing to the projective limit as $n\ra \infty$ and using 
Proposition~\ref{fin-inf} we obtain an abelian category $Coh_\gamma(X)$.
In order to get an explicit description of this category, we pick an 
affine cover $(U_i)_{i\in I}$ of $X$. All finite intersections  
$U_{i_1}\cap\dots \cap U_{i_n},\,\,\,i_1,\dots,i_n\in I,\,n\ge 1$     
are automatically affine because $X$ is separated. 
An object $\E$ of   $Coh_\gamma(X)$ is given by a
collection of  objects $(\E_i)_{i\in I}$ of categories
$\a_{U_i,\gamma}-\mod$,    together with isomorphisms of restrictions
$\phi_{ij}:(\E_i)_{|U_i\cap U_j} \simeq  (\E_j)_{|U_i\cap U_j}$ for each
pair $i,j\in I$, satisfying the associativity constraint on $U_i\cap
U_j\cap U_k$ for all $i,j,k\in I$.

The category $Coh_\gamma(X)$ seems to be in general quite degenerate. 
It has in a sense much less objects than $Coh(X)$. 
In particular, if $X$ is a two-dimensional abelian
variety and $\g$   is a non-zero translationally invariant Poisson structure,
then    vector bundles on $X$ with non-vanishing second Chern class admit no
deformations  to  objects of $Coh_\gamma(X)$. 
 
For the construction of stacks of algebroids it is more natural to 
consider as a semiclassical data not only a Poisson structure on a
manifold $X$, but a more complicated geometric  structure including  
infinitesimal deformations of $X$ (classified in the first order by 
$H^1(X,T_X)$) and also elements of $H^2(X,\O_X)$.     
The reason is that the formality morphism gives a classification 
of deformations of stacks of algebroids by gauge equivalence classes 
of solutions of the Maurer-Cartan equation in the ``derived direct image'' 
$R\Gamma(X,\wedge^\bullet T_X[1])$ (see e.g. \cite{hs}) of the
sheaf of graded Lie algebras  of polyvector fields on $X$. The locus of
the semiclassical data for which we get a sheaf of algebras (and not only a
sheaf of algebroids) is quite complicated; see \cite{nt} for calculations in
the purely symplectic case.          

\section{Affine space}

Let $X\simeq \A^n$ be an affine space, denote by $V$ the underlying 
vector space. In this section we shall discuss examples of actual and 
semi-formal quantizations of $X$  corresponding to Poisson structures on $X$
$$\g=\sum_{i,j=1}^n \g^{ij} \d_i\wedge \d_j$$
whose coefficients $\g^{ij}$ are polynomials on $V$  of low degrees. 
We shall use  canonical formulas for the $*_\h$-product in the 
standard affine structure on $X$.

\subsection{Constant coefficients}
If $\g^{ij}$ are constants then the canonical $*_\h$-product 
 is the  Moyal product
$$f *_\h g= \exp\left( \h\sum \g^{ij} \d_i\otimes \d_j\right)
(f\otimes g)_{|Diag}$$
For any given polynomials $f,g\in\O(X)=Sym(V^*)$ the product $f*_\h g$
considered  as series in $\h$ is  \emph{polynomial}, i.e. we obtain
 an actual deformation of $\O(X)$ parametrized by $\A^1=Spec(\k[\h])$.
We can evaluate the product $*_\h$ at $\h=1$. The resulting algebra
 is isomorphic to the tensor product over $\k$ of the Weyl algebra of rank
equal to $rank(\gamma)$ with a polynomial algebra. 

\subsection{Linear coefficients}
If the coefficients  $\gamma^{ij}$ are linear forms on $V$ then 
$\gl:=V^*$ is a Lie algebra and $\gamma$ is Darboux-Kirillov bracket.
 Again, as in the case of constant coefficients, the $*_\h$-product
 is a polynomial in $\h$ and can be evaluated at $\h=1$.
 The resulting algebra is isomorphic to the universal
 enveloping algebra $U(\gl)$ (see ~\cite{k1}).

\subsection{Quadratic coefficients}

Any quadratic Poisson structure is invariant under the diagonal action
of the reductive group $\gm$ on $X$  (and of the compact group $U(1)$
in the case $\k=\C$). The formulas for the canonical $*_\h$-product 
are invariant under all affine transformations. Thus, from a quadratic 
Poisson bracket we obtain a  $*_\h$-product which is homogeneous of 
degree $0$ with respect to the grading  by degree in the polynomial 
algebra $Sym(V^*)$:
$$*_\h: Sym^k(V^*)\otimes Sym^l(V^*)\otimes \k[[\h]] 
\ra Sym^{k+l}(V^*)\otimes \k[[\h]]$$

As a corollary, we conclude that the dense subspace 
$$A_\g:=Sym(V^*)\otimes \k[[\hbar]]$$  
of $Sym(V^*)[[\hbar]]$ is closed under the multiplication $*_\h$.  
The algebra $A_\g$ over $\k[[\h]]$ is a semi-formal deformation of 
$Sym(V^*)$, free as a $\k[[\h]]$-module.  Inverting $\hbar$ we obtain a
finitely generated algebra  $$A'_\g:=A_\g\otimes_{\k[[\hbar]]} \k((\hbar))$$
over the field $\k((\hbar))$ satisfying an analog  of the
Poincar\'e-Birkhoff-Witt theorem. It follows from results of Drinfeld
(see \cite{d2}) that the algebra $A'_\g$ is a Koszul quadratic algebra,  
noetherian and of cohomological dimension $n$.

Restriction of the $*_\h$-product to the case $k=l=1$ gives a 
$\k[[\h]]$-linear map 
$$V^*\otimes V^*\otimes \k[[\h]]\ra Sym^2(V^*)\otimes \k[[\h]]$$
This is an epimorphism (because it gives the obvious epimorphism
at $\h=0$). Let us denote by $R_\h^{(2)}$ its kernel.
It is easy to prove that $A_\g$ is generated by free $\k[[\h]]$-module 
$V^*\otimes \k[[\h]]$ with quadratic relations $R_\h^{(2)}$. The module 
$R_h^{(2)}$ is projective of rank equal to $\dim(\wedge^2(V))$.
Analogously, we define the module $R_\h^{(3)}$ as the kernel of 
the triple $*_\h$-product 
$$V^*\otimes V^*\otimes V^*\otimes\k[[\h]]\ra Sym^3(V^*)\otimes \k[[\h]]$$
Define a $GL(V)$-invariant scheme $QAlg_V\subset Gr(a_2,V)\times Gr(a_3,V)$
where  $$a_k:=\dim(V^{\otimes k})- \dim(Sym^k(V)),\,\,\,k=2,3\,,$$ as the
closed  subscheme consisting of pairs $\left(R^{(2)},R^{(3)}\right)$ such
that  $$R^{(3)}\supset R^{(2)}\otimes V + V\otimes R^{(2)}$$
We have a point $Comm_V\in QAlg_V$ corresponding to the standard 
commutative algebra $Sym(V^*)$. Deformation quantization gives a formal path 
in  $QAlg_V$ starting at $Comm_V$. 

The formal completion $\widehat{QAlg}_V$ of  scheme $QAlg_V$ at point 
$Comm_V$ can be considered as the formal moduli space of quadratic algebras.  
It follows from results of Drinfeld (see \cite{d2}) that $\widehat{QAlg}_V$ 
is the universal base of a formal flat deformation of $Sym(V^*)$ considered 
as a graded associative algebra with a \emph{given} space $V^*$ of generators. 

Analogously, we define $GL(V)$-invariant scheme $QPois_V\subset
\wedge^2(V)\otimes Sym^2(V^*)$  as the closed subscheme of quadratic bivector
fields on $V$ giving a Poisson  structure. This is obviously a \emph{cone}
because the defining equation 
$$[\gamma,\gamma]=0\in \wedge^3(V)\otimes Sym^3(V^*)$$
is homogeneous. We denote by $\widehat{QPois}_V$ the formal completion 
of $QPois_V$ at $0$. Deformation quantization gives a formal
$GL(V)$-equivariant map
$${quant}_V:\widehat{QPois}_V\ra \widehat{QAlg}_V\,\,,$$
the formal path in $QAlg_V$ corresponding to $*_\h$ is the image of the formal
completion at $\h=0$ of the straight path  $\h\mapsto \gamma\times\h$ in the
cone $QPois_V$. 

\begin{prop} The map ${quant}_V$ is an isomorphism.
\end{prop}
 
{\it Proof:} A point (over an Artin algebra) of $QAlg_V$ infinitesimally close
to $Comm_V$  gives a deformed product on $Sym(V^*)$. In order to define this
product we notice that the deformed algebra is generated by $V^*$, thus it is 
a quotient of the tensor algebra. Using the standard embedding of subspaces
$Sym(V^*)\hookrightarrow \otimes^\bullet(V^*)$  we obtain the identification
of the symmetric algebra with the deformed one, thus we get a deformed
product on $Sym(V^*)$. Applying  arbitrary $\gm$-equivariant left inverse to 
the $L_\infty$-quasi-isomorphism\footnote{this $L_\infty$-morphism is the  
composition of the formality $L_\infty$-morphism 
$T_{poly}^\bullet(X)\ra D_{poly}^\bullet(X)$ with the obvious embedding 
of $D_{poly}^\bullet(X)$ in the cohomological Hochschild complex,
see \cite{k1} for the notations.}
$$T_{poly}^\bullet(X)\ra C^\bullet(\O(X),\O(X))[1]$$
we obtain a formal ``dequantization'' map  $${dequant}_V:
\widehat{Qalg}_V\ra \widehat{QPois}_V$$ Both maps  ${quant}_V$ and
${dequant}_V$ induce (mutually inverse) isomorphisms of the tangent
spaces at base points.

\begin{lemma} Let $Y/\k$ be a formal affine scheme with one closed point 
$y$ and $f:Y\ra Y$ be a map inducing the identity map on $T_y Y$. 
Then $f$  is invertible.
\end{lemma}

 {\it Proof of the Lemma:} Let us denote by $J$ the maximal ideal 
in $\O(Y)$ corresponding to closed point $y$. The map $f^*$ preserves
complete decreasing filtration on $\O(Y)$ by powers of $J$.  The induced map
on the associated graded algebra is the identity because the quotient spaces
$J^n/J^{n+1}$ are identified with quotient spaces of the $n$-th symmetric
power of the cotangent space $T^*_y Y$.  Therefore $f^*$ is an isomorphism.
\qed

Applying the lemma to maps ${quant}_V\circ {dequant}_V$  and 
${dequant}_V\circ {quant}_V$ we obtain the invertibility of ${quant}_V$. \qed
 
It is not clear how the map ${quant}_V$ looks like even in the simplest
case $\dim(V)=2$. For example, a natural guess is that  the canonical
quantization of the Poisson bracket 
$$\hbar \cdot xy \,{\partial\over \partial  x} \wedge {\partial\over \partial 
y}$$  is quadratic algebra given by relations 
$$y  *_\hbar x = exp(\hbar) xy$$
but it is not clear how to prove it.

\begin{conj} For the case $\k=\C$ and for the canonical deformation
quantization as in \cite{k1} or \cite{t}, the series for the map $quant_V$ 
is convergent in   $\C$-topology in a neighborhood of zero. 
\end{conj}

It is well-known that for Sklyanin brackets the quantized product is
parametrized  by a point on an elliptic curve.  We expect that in this case 
the map $quant_V$ is inverse to an elliptic function.

\begin{quess} Can one in general express the inverse to the map $quant_V$
 in terms of periods (Gauss-Manin connection)?
 Is the map $quant_V$ independent of the choice of constants in universal
deformation quantization?
 If so, is the series $quant_V$ convergent in a disc in the p-adic case,  
and is it an entire function in the case $\k=\C$? 
\end{quess}

\subsection{Nonhomogeneous quadratic coefficients} 

The definition of the quantized algebra (free as a $\k[[\h]]$-module)
can be repeated without difficulties  for a general Poisson bracket $\g$ 
whose coefficients are polynomials of degree $\le 2$.
The resulting algebra is no longer graded, it is filtered by  free
 $\k[[\h]]$-modules of finite rank. As above, we obtain a semi-formal
deformation of $\O(X)$ endowed with a natural filtration. The associated
graded algebra  coincides with the canonical semi-formal quantization of the
purely quadratic part of $\g$.

Moreover, one can reduce the study of quantizations of nonhomogeneous
quadratic  Poisson brackets to the case of purely quadratic brackets by adding 
one central variable.

\subsection{Cubic brackets etc.}  
Some Poisson brackets of higher degree can be quantized semi-formally 
and one obtain algebras over $\k[[\h]]$, free as $\k[[\h]]$-modules. 
This can be done, for example, for the following bracket on $X=\A^3$: 
$$\{x_1,x_2\}=P(x_3),\,\,\,\{x_1,x_3\}=\{x_2,x_3\}=0$$
where $P\in \k[x]$ is an arbitrary polynomial.
However, it seems that  there is no way to construct semi-formal 
quantizations for  generic Poisson brackets on  affine spaces, including
terms of degree $3$ and higher. In Section 4 we propose an explanation of
the impossibility of the quantization. 

\section{Projective manifolds}

\subsection{Quantum homogeneous coordinate ring} 
 
Let $(X,\g)$ be a smooth projective variety with Poisson structure, and $\L$ 
be an ample line bundle on $X$. The question is whether one can quantize 
$X$ ``preserving'' the bundle $\L$. More precisely, we would like to 
quantize the homogeneous coordinate ring 
$$A^\bullet:=\sum_{n\ge 0}\Gamma(X,\L^{\otimes n})$$ 
Here (because of the compactness of $X$)  there is no difference between
formal and semi-formal quantizations of the graded algebra $A^\bullet$, its
graded components are already finite-dimensional, and also $A^\bullet$ is
finitely generated.  We shall see that in fact it is possible to deform only a
subalgebra of $A^\bullet$ obtained from  $A^\bullet$ after the removal of a
finite number of homogeneous components.

\begin{dfn} A line bundle $\L$ is \emph{formally quantizable}  
iff there  exists a Poisson structure $\tg$ on the total space $\L^\times$  
of the principal $\gm$-bundle associated with $\L$, such that  
$\tg$ is $\gm$-invariant and projection $\pi:\L^\times\ra X$ is a Poisson map.
\end{dfn}

\begin{const} Assume that   $H^1(X,\O_X)=H^2(X,\O_X)=0$.
 Then for any Poisson structure $\gamma$, given a formally quantizable 
ample line bundle $\L$ and a choice of $\tg$, we construct a canonical 
noetherian graded algebra $A^\bullet_{\widetilde{\g}}$ flat over $\k[[\h]]$ 
(and also free as $\k[[\h]]$-module), deforming the graded algebra 
$A_{(tr)}^\bullet:=\k\cdot 1\oplus \bigoplus_{n\ge n_0}\Gamma(X,\L^{\otimes
n})$ where $n_0\ge 1$ is sufficiently large.\footnote{The subscript $(tr)$ in
the notation indicates that $A_{(tr)}^\bullet$ is a {\em truncated} version of
the  homogeneous coordinate ring of $X$.} 
\end{const}  

Let us cover $X$ by  affine open subschemes $(U_i)_{i\in I}$.
For every $i\in I$ there exists a $\gm$-invariant connection $\nabla_i$ 
on the tangent bundle restricted to 
$\widetilde{U}_i:=\pi^{-1}(U_i)\subset \L^\times$,
because group scheme $\gm$ is reductive. Let us choose  such 
connections $\nabla_i$ for all $i\in I$.   Thus, we get a $\gm$-invariant
product $*_{\hbar,i}$ on $\O(\widetilde{U}_i)[[\h]]$ for any $i\in I$, an
isomorphism of algebras $$g_{ij}:( \O(\widetilde{U}_i\cap
\widetilde{U}_i)[[\h]],*_{\hbar,i}) \ra  ( \O(\widetilde{U}_i\cap
\widetilde{U}_i)[[\h]],*_{\hbar,j}) $$ 
 given by an infinite series in $\h$ of $\gm$-invariant differential
operators,  and certain invertible elements (see Appendix, part A.1) 
$$a_{ijk}\in \left((\O(\widetilde{U}_i\cap \widetilde{U}_j\cap
\widetilde{U}_k)[[\h]])^{\times}\right)^{\gm}\simeq
(\O({U}_i\cap {U}_j\cap {U}_k)[[\h]])^{\times}$$ 
satisfying a twisted multiplicative 2-cocycle condition. 
At $\hbar =0$ the products $*_{\hbar,i}$ coincide with the standard product, 
isomorphisms $g_{ij}$ are identity maps and all elements 
$a_{ijk}$ are equal to $1$.

 We claim that  after a $\gm$-equivariant
gauge transformation (modifying isomorphisms $g_{ij}$) one can make all
$a_{ijk}$ equal to $1$.   The proof goes term by term, the lowest (in
powers of $\h$)  non-zero term in $\log(a_{ijk})$ is a \Cech $3$-cocycle of 
the covering  $(U_i)_{i\in I}$ with coefficients in $\O_X$, we kill 
it using the assumption $H^2(X,\O_X)=0$.

 After making all elements $a_{ijk}$ equal to $1$ we obtain 
$\gm$-equivariant gluing data for associative 
algebras $\O(\widetilde{U}_i)[[\h]]$ (with products $*_{\hbar,i}$). 
Thus, we can associate with it a  \Cech complex
$$\sum_{i\in I} \O(\widetilde{U}_i)[[\h]]\ra 
\sum_{i,j\in I} \O(\widetilde{U}_i\cap \widetilde{U}_i)[[\h]]\ra \dots$$
The differential in this complex is $\gm$-equivariant and $\k[[\h]]$-linear, a
deformation  of the \Cech complex of the covering $(\widetilde{U}_i)_{i\in I}$
 of $\L^\times $ with coefficients in $\O_{\L^\times}$.
 It follows from the ampleness of $\L$  that for some sufficiently large
$n_0$, all $n\ge n_0$ and all $i\ge 1$ we have  $H^i(X,\L^{\otimes n})=0$. Then
by the semicontinuity of cohomology we conclude  that the \Cech complex
associated with the gluing data for homogeneous components (with respect to
$\gm$-action) of   $\O(\widetilde{U}_i)[[\h]]$ of degrees $\ge n_0$ has
cohomology only in the cohomological degree $0$. Denote by 
$A^n_{\widetilde{\g}}$ the $0$-th cohomology of the $n$-th graded component of
the \Cech complex from above.   It follows from the invariance of the Euler
characteristic under deformations  that for $n\ge n_0$, the space
$A^n_{\widetilde{\g}}$ is a  free finitely generated module over $\k[[\h]]$. 
The graded module over $\k[[\h]]$ 
$$A^\bullet_{\widetilde{\g}}:=\k[[\h]]\cdot 1\oplus \bigoplus_{n\ge n_0} 
A^n_{\widetilde{\g}}$$  forms
obviously an associative algebra, a flat deformation over $\k[[\h]]$ of the 
truncated graded algebra $$A_{(tr)}^\bullet=\k\cdot 1 \oplus \bigoplus_{n\ge
n_0}\Gamma(X,\L^{\otimes n})$$

 Analogously, using the assumption $H^1(X,\O_X)=0$ one can show that for
different choices of ways to kill  $\log(a_{ijk})$, and 
also for different choices  of connections on 
$\widetilde{U}_i$  we get isomorphic algebras, with automorphisms defined
up to an inner automorphism.  The invertible element involved in the
 above-mentioned inner automorphism is automatically an element of
 $A^0_{\widetilde{\g}}$. In the non-deformed case we get only  constants, 
$H^0(X,\O_X)=\k\cdot 1$. Thus, by semicontinuity after the deformation,
we can get only a smaller space. On the other hand,
constants should remain central elements after the deformation. 
Thus we see that the ambiguity in the definition of isomorphisms between 
variously constructed graded associative algebras is absent.

Finally, the noetherian property of $A^\bullet_{\widetilde{\g}}$ follows 
easily from Proposition~\ref{Noether}.

Notice that our above construction can be interpreted  as a construction of a 
semi-formal deformation of a {\em singular} Poisson scheme
$(Spec(A_{(tr)}^\bullet),\widetilde{\g})$.

\subsection{Quantum coherent sheaves: another approach}

In the commutative case any line bundle $\L$ on $X$ gives an autoequivalence 
 $F:=\L\otimes \cdot$ of the abelian category $Coh(X)$.  
Also,  the category $Coh(X)$ contains a distinguished object, the structure 
sheaf $\O_X$. The pair $(F,\O_X)$ gives rise to a $\Z$-graded algebra
$$A^\bullet=\oplus_{n\in \Z} A^n,\,\,\,A^n:=Hom(\O_X, F^n(\O_X))=
\Gamma(X,\L^{\otimes n}).$$
By Serre's theorem, for ample $\L$ the category $Coh(X)$ is naturally 
equivalent to the  quotient category  of the category of finitely 
generated graded $A^\bullet$-modules by the subcategory of
finite-dimensional modules.  It is easy to see that one can replace algebra
$A^\bullet$ by a slightly  truncated version $A^\bullet_{(tr)}$ as above and
obtain the same category $Coh(X)$.

 One of approaches to noncommutative algebraic geometry (see e.g.~\cite{az})
 exploits Serre's theorem using it as a seed for the {\em definition} of a
noncommutative projective variety     as an abelian category endowed with an
autoequivalence and a distinguished  object satisfying  some natural
conditions. So, it seems reasonable  to define the quantized category of
coherent sheaves $Coh_{\widetilde{\g}}(X)$ as the quotient of the category
of  finitely generated graded $A^\bullet_{\widetilde{\g}}$-modules  
by the subcategory of finite-dimensional modules. It has a distinguished 
object $\O_{X,\widetilde{\g}}$ and an autoequivalence 
$(\L\otimes\cdot)_{\widetilde{\g}}$.
  
Notice that in Section 1.5 we already defined an abelian category
  $Coh_{\g}(X)$ using gluing of finitely generated modules.
 In the case of the trivial Poisson structures $\g,\widetilde{\g}$ both
categories  $Coh_{\widetilde{\g}}(X)$ and $Coh_{\g}(X)$ are naturally
equivalent  to the category of coherent sheaves on the scheme
$X\times_{Spec\,(\k)}  Spec\,(\k[[\h]])$. This fact motivates the
following question:

\begin{quess} Are categories $Coh_{\widetilde{\g}}(X)$ and  
$Coh_{\g}(X)$ naturally equivalent?
Is the cohomological dimension of category $Coh_{\widetilde{\g}}(X)$ finite 
(in fact equal to $\dim(X)$)? 
\end{quess}

\subsection{Dependence on the choice of the lift $\widetilde{\g}$} 

Let $\widetilde{\g}_1$ and $\widetilde{\g}_2$ be two $\gm$-equivariant
 Poisson structures on $\L^\times$ extending the same Poisson structure $\g$
on $X$. We want to compare graded algebras $A_{\widetilde{\g}_1}^\bullet$ and
 $A_{\widetilde{\g}_2}^\bullet$, and also corresponding quantized
 categories of coherent sheaves. It is easy to see that the difference 
$\widetilde{\g}_1-\widetilde{\g}_2$ can be (informally) written as $e\wedge v$
 where $e$ is the infinitesimal generator of $\gm$-action on $\L^\times$
 and $v$ is a vector field on $X$ preserving $\g$.
 Thus, we see that on the semiclassical level the algebra
$A_{\widetilde{\g}_2}^\bullet$  is a twisted version of
$A_{\widetilde{\g}_1}^\bullet$ associated with an automorphism of the graded
algebra $A_{\widetilde{\g}_1}^\bullet$ (see \cite{av} for a discussion of
twisted graded algebras in the commutative case).  It is natural to expect that
the same  holds after quantization. Translation of this into  the categorical
language means that one can (conjecturally) canonically identify abelian
categories $Coh_{\widetilde{\g}_1}(X)$ and  $Coh_{\widetilde{\g}_2}(X)$
together with distinguished objects    $\O_{X,\widetilde{\g}_1}$ and
$\O_{X,\widetilde{\g}_2}$.  
Endofunctors $(\L\otimes\cdot)_{\widetilde{\g}_1}$ and 
$(\L\otimes\cdot)_{\widetilde{\g}_2 }$ become different after this
identification.

\subsection{Quantum projective spaces}

Let $X$ be a projective space $\P(V)$ where $V$ is a finite-dimensional vector
space, $\dim(V) \ge 2$, $\g$ be a Poisson structure on $X$ and $\L$ be an
ample line bundle on $X$ (i.e.  $\L\simeq\O(m)$ for some $m>0$). First of all,
it is easy to see that if $\L$ is formally quantizable then  also $\O(1)$ is
 formally quantizable and sets of lifts $\widetilde{\gamma}$ in both bundles 
can be canonically identified. Moreover, corresponding graded algebras are 
related in the obvious way, one is obtained from another  by keeping all 
graded components with degrees divisible by $m$.  Thus, we can assume safely 
that $m=1$ and $\L=\O(1)$. 

Any $\gm$-equivariant Poisson structure $\widetilde{\g}$ on
$\L^\times=V\setminus\{0\}$ extends by Hartogs
principle\footnote{ $\O(X\setminus Z)=\O(X)$ for any smooth $X$ and closed
 $Z\subset X$ of codimension $\ge 2$} to a
quadratic  Poisson structure on $V$.  Here we have an alternative construction
of the quantized graded algebra  (see Sect. 2.3). It is not hard to identify
canonically the outputs of both  constructions.

Finally, for any Poisson structure on $\P(V)$, the line bundle $\O(1)$ is 
formally quantizable (A.~Bondal, private communication).  
Thus, for projective spaces quantizations are governed by quadratic algebras.

\subsection{Non-quantizable manifolds and vanishing conditions} 

Let $X$ be a projective symplectic manifold 
 with the obvious non-degenerate Poisson structure $\g$.
 We claim that $(X,\g)$ can not be 
quantized by the general procedure described above in 3.1. There are two
independent reasons. First of all, any formally quantizable  line bundle on $X$
 carries a flat connection,  thus it cannot be ample. 
Secondly, the vanishing condition $H^1(X,\O_X)=H^2(X,\O_X)=0$ fails.  We
consider the first reason to be much more serious. It means that one cannot
 obtain a deformation of the homogeneous coordinate ring
$A^\bullet$ (or even  $A^\bullet_{(tr)}$).

The condition $H^1(X,\O_X)=H^2(X,\O_X)=0$ in general seems to be almost 
redundant. The existence  of the Poisson structure already gives a
restriction to the existence of  a nonvanishing global $2$-form. Thus, we
think that our conditions  for the semi-formal quantizability  of projective
Poisson manifolds  (i.e. for the existence of a (semi-)formal deformation of
the truncated homogeneous  coordinate ring) are not far
from the best possible.

In the Introduction we mentioned that it is impossible to quantize in a proper
sense (i.e. semi-formally) abelian varieties.  The reason is that there is no
formally quantizable  ample line bundle on an abelian variety $X$ endowed with
any non-zero Poisson structure. 

  \section{Affine manifolds}

 We start with a technical notion which will be used later.
 Let $(Y,\gamma)$ be an algebraic manifold with bivector field $\gamma\in
\Gamma(Y,\wedge^2 T_Y)$, and $Z\subset Y$ be a reduced divisor. 

\begin{prop}
The following three conditions are equivalent: 
\begin{itemize}
\item the sheaf of ideals $I_Z$ of functions vanishing at $Z$ is also an ideal
with respect to the  bracket $\{f,g\}:=\langle df\otimes dg,\gamma\rangle$,
\item the same is true for a sheaf $J$ of ideals such that $rad(J)=I_Z$,
\item there exists a Zariski dense subset $U\subset Z$ consisting of smooth
points of $Z$ such that for any $p\in U$ and local coordinates $(x_1,\dots
x_n)$
 near $p$ such that $Z$ is given by equation $x_1=0$, bivector field $\gamma$ 
belongs to $\O_Y$-submodule generated by exterior products of
 fields $x_1 \partial/\partial x_1$ and $\partial/\partial x_i,\,\,\, i\ge 2$.
\end{itemize}
\end{prop}

This proposition follows easily from the Hartogs principle.
We say that $\gamma$ is {\em tangent} to $Z$ iff one of the three equivalent
conditions above is verified.  

Let $(X,\gamma)$ be an affine Poisson manifold.

\begin{dfn}  A smooth compactification ${\overline X}\supset X$
 is {\em quantizable} iff  
\begin{itemize}
\item $H^1({\overline X},\O_X)=H^2({\overline X},\O_X)=0$,
\item the bivector field $\gamma$ extends to ${\overline X}$,
\item $\gamma$ is tangent to divisor $Z:={\overline X}\setminus X$,
\item the class of divisor $Z$ in $Pic({\overline X})$ is ample.
\end{itemize}
\end{dfn}

The first condition does not depend on the compactifications, by the birational
invariance of Hodge numbers $h^{0,\bullet}$. It holds e.g. for rational
varieties. A quantizable compactification gives an increasing filtration by
finite-dimensional spaces on $\O(X)$:
$$\O(X)=\cup_{m}\O(X)_{\le m},\,\,\,\,\,\O(X)_{\le m}:=\Gamma({\overline
X},\O(m Z))$$ This filtration is compatible with the product and the Poisson
structure in the following sense:
  $$\O(X)_{\le m_1}\cdot \O(X)_{\le m_2}\subset \O(X)_{\le
m_1+m_2}\,,$$ $$\{\O(X)_{\le m_1}, \O(X)_{\le m_2}\}\subset \O(X)_{\le
m_1+m_2}$$  Conversely, if such a filtration is given and it is generated in
the multiplicative sense by its finite part:
$$\exists m\ge 1, \,\forall N\ge 1\,\,\,\,\O(X)_{\le N}=
\sum_{i_1,\dots,i_k\le m,\,\,i_1+\dots +i_k=N}
\O(X)_{\le i_1}\cdot \dots \cdot\O(X)_{\le i_k}
$$
 then it gives a compactification of $X$ by a normal variety, with the 
Poisson structure tangent to the smooth part of the compactification divisor. 
Thus, we get a compactification which is quantizable iff it is smooth.
Notice that our condition is considerably weaker than the usual one 
$$\{\O(X)_{\le m_1}, \O(X)_{\le m_2}\}\subset \O(X)_{\le m_1+m_2-1}$$ which
 holds for semiclassical approximations to algebras of differential operators
and to universal enveloping algebras.

 \begin{const} With any quantizable compactification ${\overline X}$ of 
a Poisson manifold $(X,\gamma)$ we associate a canonical algebra 
$\O_{\gamma,{\overline X}}(X)$ which is 
a semi-formal deformation of $\O(X)$ over $\k[[\hbar]]$.
 \end{const}

Let us denote by $\L$ the ample line bundle associated with divisor $Z$.
 By its definition, the bundle $\L$ is endowed with a section $s$ vanishing 
exactly at $Z$. Moreover, the trivial lift of the Poisson structure 
$\gamma$ on $X$ to the pullback of $X$ extends without poles to the whole 
total space $\L^{\times}$. Applying the construction of the quantization in 
the projective case we obtain a graded algebra 
$A^\bullet_{\widetilde{\g}}$ over $\k[[\h]]$.

 We claim that there exists a canonical $\gm$-equivariant embedding of a 
truncated polynomial algebra  
$$\k[[\h]]\cdot 1\oplus \bigoplus_{n\ge n_0} \k[[\h]]\cdot s^n$$
into the center of $A^\bullet_{\widetilde{\g}}$. At $\h=0$   it is given  
by non-negative powers of the section $s$. The reason of the existence of 
such an embedding is that in deformation quantization one gets a 
canonical identification (preserving the cup-product) between  the 
Poisson cohomology of the Poisson manifold and the Hochschild cohomology  
of the quantized algebra (see~\cite{k1}, Section 8.2). 
In particular, for $0$-th cohomology one obtains an algebra isomorphism 
between centers  in the Poisson and in the associative sense. 
Strictly speaking, all what was said is true in the $C^{\infty}$-case, but
one can construct easily analogous  isomorphisms for non-affine algebraic
varieties, with Poisson cohomology replaced by hypercohomology, and a similar
replacement for the Hochschild cohomology. The global $0$-th Poisson
cohomology of $\L^{\times}$  contains the subalgebra $\oplus_{n\ge 0}\k\cdot
s^n $.  Thus, the same should hold for $A^\bullet_{\widetilde{\g}}$. 

The conclusion is that the algebra $A^\bullet_{\widetilde{\g}}$ contains 
a multiplicative submonoid 
$$M:=\{s^n|\,n=0\mbox{ or } n\ge n_0\}$$ consisting of homogeneous (with
respect to $\Z$-grading) central elements. We define quantized algebra 
$\O_{\gamma,{\overline X}}(X)$ as the homogeneity degree $0$ part of the 
localization $A^\bullet_{\widetilde{\g}}[M^{-1}]$ of 
$A^\bullet_{\widetilde{\g}}$ with respect to $M$. This quantized algebra 
is filtered by finitely generated free $\k[[\h]]$-modules.
 
Notice the important difference between the rigidifying data here (a 
quantizable compactification) and in section 1.2 (a connection on the tangent
bundle). 

\subsection{Examples and questions}

First of all, for a polynomial Poisson structure $\g$ on $\A^n$ the standard
 compactification by $\P^n$ is quantizable iff the coefficients of $\g$ are
polynomials of degree $\le 2$. The quantization $(\A^n,\g)$
constructed using the procedure proposed in this  section gives the same
algebra as the one described in the section devoted to the quantization of the
affine space. For general Poisson structure on $\A^n$ with coefficients
including terms of higher degree, it is very unlikely to have a quantizable
compactification. In general, bivector fields tend to acquire
poles after blow ups.

One can try also to quantize a Poisson submanifold of the affine space
endowed with a bracket of degree $\le 2$. For example, a natural guess 
 would be that the quantization of a
regular coadjoint orbit  of a semisimple algebraic group $G$ coincides with the
appropriate quotient  of $U(Lie(G))$ given by the Harish Chandra isomorphism.

Another source of algebraic Poisson manifolds is provided by theory of
``algebraic'' quantum groups,  see e.g. \cite{so}. For example, 
 let us realize the standard Poisson-Lie group $SL_2$ as a quadric hypersurface
in the space of matrices $\A^4$, and then compactify $SL_2$ by its closure 
in $\P^4$. We expect that the corresponding quantized algebra coincides with
the algebra of  regular functions on quantum $SL_2$.

It seems that the criteria proposed here for semi-formal quantizability of 
affine Poisson manifolds, are  too strong. 
For example, if the compactification has singularities in codimension $4$, 
our construction can still be performed.

\begin{quess} What are natural necessary and sufficient conditions
 for quantizability of smooth Poisson varieties which are neither
 projective nor affine (like nilpotent coadjoint orbits or
Poisson-Lie homogeneous spaces $G/N$)? Also, what can one say about 
the quantization of spaces with singularities?
\end{quess}

It is known after results of J.~Dixmier, H.~W.~E.~Jung and
L.~Makar-Limanov  (see \cite{bw} and references therein) that in the case 
$X=\A^2$ (with coordinates $x,y$) and 
$\g={\partial\over\partial x} \wedge {\partial\over\partial x}$ 
the group of automorphisms of $(X,\g)$  coincides
with the group of automorphisms of the Weyl algebra $W_1$. Thus,
in this case the quantized algebra can be constructed {\em canonically}. There
is another example (see \cite{rsz}) of a huge group of automorphisms  of an 
algebraic Poisson variety acting on a (candidate to) quantization. In this
example the Poisson variety (it has singularities) is the  variety of 
representations in $SL_2$ of the fundamental group of a surface.  The group
of symmetries contains the mapping class group, and the quantization 
constructed in \cite{rsz} is equivariant.

\begin{ques} Is it possible to identify canonically the algebras 
$\O_{\gamma,{\overline X}}(X)$ for all quantizable
compactifications of $(X,\g)$? 
\end{ques}
 
A possible approach to the positive solution would be to intersect
 filtrations on $\O(X)$ corresponding to two different quantizable
compactifications and try to use the obtained filtration in order to make an
``interpolation'' between two quantizations.  An obstacle in this approach is
that the compactification arising from the intersected  filtration is not
necessarily smooth.

A positive answer will be probably interesting from the point of view of
 the orbit method. 

\appendix

\section{Algebroids in deformation quantization}

\subsection{Combinatorial description  of algebroids}

Suppose that we are given a simplicial complex $X$ 
 with the set of vertices $I$, an identification $i\mapsto \E_i$
between $I$  and the set of object of an algebroid $\a$, and an isomorphism
 $f_{ij}\in Hom_\a(\E_i,\E_j),\,\,\,f_{ji}=f_{ij}^{-1}$ for every
 $1$-dimensional face $\{i,j\}\subset I$ of $X$. 
 Then we obtain the following data:
\begin{itemize}
\item for every $i\in I$ an associative unital
algebra  $\a_i:=End_\a(\E_i)$, 
\item for every pair $(i,j)$ of vertices such that $\{i,j\}$ is an edge 
of $X$, an isomorphism of algebras $g_{ij}:\a_i\ra \a_j$, 
$g_{ji}=g_{ij}^{-1}$,
\item for every triple $(i,j,k)$ such that $\{i,j,k\}$ is a 
$2$-dimensional face of $X$, an invertible element $a_{ijk}\in \a_i$.
\end{itemize}
Namely, the isomorphism $g_{ij}$ is given by 
$g_{ij}(h):=f_{ij}\circ h\circ f_{ji}$, and the element 
$a_{ijk}$ is equal to $f_{ki}\circ f_{jk}\circ f_{ij}$. 
This data satisfies the following constraints:
 $$a_{ikj}=a_{ijk}^{-1},\,\, a_{jki}=g_{ij}(a_{ijk}),\,\,\,
g_{ki}\circ g_{jk}\circ g_{ij}=Ad(a_{ijk}^{-1})$$
and the following ``tetrahedron'' equation\footnote{The data and conditions 
presented above are familiar from Giraud theory of gerbs and non-abelian 
cohomology (see \cite{gi}), and similar to  formulas written recently by
physicists  for D-branes in the presence of B-field (see e.g. \cite{s}).}
for every $i,j,k,l$ such that  $\{i,j,k,l\}$ is a $3$-dimensional face of $X$:
$$ a_{ikl}\circ a_{ijk}=a_{ijl}\circ g_{jl} (a_{jkl})$$

Conversely, for any $3$-connected simplicial complex $X$ endowed
with  data $(\a_i,g_{ij},a_{ijk})$ satisfying conditions as above one can  
canonically reconstruct $\a$. Namely, every simplicial  path
$p=(i_0,\dots,i_n)$ gives an isomorphism  
$g_{i_{n-1} i_n}\circ \dots\circ g_{i_0 i_1}$ between  the algebras  
$\a_{i_0}$ and   $\a_{i_n}$. We denote by $\B_p$ the diagonal bimodule  
(left over $\a_{i_0}$ and  right over $\a_{i_n}$)
using this isomorphism.  Then it is easy to see that using elements $a_{ijk}$
one can construct  an identification between the bimodules $\B_p$ and 
$\B_{p'}$ for two paths  with the same endpoints, connected by a 
simplicial homotopy. The tetrahedron  condition implies that the 
identification does not depend on the choice  of homotopy. 
For any two elements $i,j\in I$ the result of the identification of  
the bimodules $\B_p$ over all  paths with ends $i,j$ is 
canonically isomorphic to the space of morphisms  $Hom_\a(\E_i,\E_j)$ in the
algebroid $\a$. 

\subsection{Families of algebras with connections} 
Here we describe an infinitesimal version of the construction of the 
combinatorial description of algebroids given above. 
In what follows manifolds are understood in a broad sense, they are
either $C^\infty$-manifolds (maybe infinite-dimensional)
or formal manifolds over a field of characteristic zero.

Let $A$ be a  vector space with  a marked vector $\1_A\in A$, and $M$ be a
contractible manifold (in fact only $3$-connectedness is necessary).  
Consider the tensor product  $${\bf g}:=\G(M,\Omega^\bullet_M)\otimes
C^\bullet(A,A)[1]$$ of the de Rham algebra of $M$ and the shifted Hochschild
complex of $A$  (endowed with zero multiplication) 
 as a differential graded Lie algebra. Suppose that we are given 
 a solution $\gamma\in {\bf g}^2$  of the Maurer-Cartan equation 
$$d\gamma+{[\gamma,\gamma]\over 2}=0$$
satisfying some additional unitality constraints (see below).
We associate with $\gamma$ an algebroid $\a$ whose objects are 
points of $M$. 
 First of all, decompose $\gamma$ into the sum of bi-homogeneous components:
$$\gamma=\gamma_0+\gamma_1+\gamma_2,$$
where $\gamma_i\in \G(M,\Omega^i_M)\otimes C^{2-i}(A,A)$ for $i=0,1,2$.
 The Maurer-Cartan equation gives $4$ constraints
\begin{enumerate} 
\item $[\gamma_0,\gamma_0]=0$,
\item $[\gamma_1,\gamma_0]+d\gamma_0=0$,
\item $[ \gamma_2,\gamma_0]+\left([\gamma_1,\gamma_1]/2+d\gamma_1\right)=0$,
\item $ [\gamma_2,\gamma_1]+d\gamma_2=0\,.$
\end{enumerate}

The unitality constraints  mentioned above are 
\begin {enumerate} 
\setcounter{enumi}{4}
\item ${\gamma_0}_{|x}(\1_A,f)=f$ for any $f\in A,\,\,\, x\in M$, 
\item $\gamma_1(\1_A)=0\in \Gamma(M, \Omega^1_M)\otimes A\,.$
\end{enumerate} 

The first and the fifth constraint say that $\gamma_0$ gives a 
family of associative products on $A$, parametrized by $M$, with unit $\1_A$.
 We define the endomorphism algebras of objects of $\a$ as
 $$\a_x:=(A,{\gamma_0}_{|x}),\,\,\,\forall\, x\in M$$
Let us introduce the connection $\nabla:=d+\gamma_1$ on the
trivial  bundle on $M$ with fiber equal to $A$.
The second and the sixth constraint mean that the holonomy
 of this connection along any path\footnote{Here we assume that the notion of
the holonomy makes sense. In the application below  the manifold $M$ and
the space $A$ are both infinite-dimensional.} preserves the structure  of an 
associative unital algebra.

Let us consider a piecewise-smooth map $\phi$ of two-dimensional disc $D^2$ to
$M$ and pick a point $p\in \partial D^2$.  Denote by $x$ the point 
$\phi(p)\in M$. Using the third constraint one can construct an element 
$a_{\phi,p}\in \a_x$ such that the monodromy along $\phi(\partial D^2)$  
is an inner automorphism of algebra $\a_x$ given by conjugation  with 
$a_{\phi,p}$. We leave to the reader as an exercise to write an explicit  
formula for $a_{\phi,p}$.
[Hint:   decompose $D^2$ into small triangles, use  ideas from the 
combinatorial version to define an approximate formula for $a_\phi$, 
then go to the limit making triangles smaller and smaller.]

 Also, it is easy that the fourth constraint implies
 that $a_{\phi,p}$ depends only on the loop $\phi_{|\partial D^2}$ and the
 tetrahedron equation holds.  
 
Now we define the space of morphisms $Hom_\a(x,y)$ for an arbitrary pair
of  points $x,y\in M$. 
  For any  path $\psi:[0,1]\ra M,\,\,\,\psi(0)=x,\,\psi(1)=y$ we
 identify the algebras $\a_x$ and $\a_y$ using the holonomy of $\nabla$ 
along $\psi$. THe space $Hom_\a(x,y)$ will be identified with the diagonal 
bimodule with $\a_x$. For two different paths  $\psi_1,\psi_2$ the 
identifications of $Hom_\a(x,y)$ differ by the right multiplication by 
$a_{\psi_1^{-1}\circ\psi_2, 0}$. The tetrahedron identity ensures that the 
space $Hom_\a(x,y)$ is well-defined. The composition of morphisms is defined 
in an obvious way using the holonomy along paths.
 
\begin{rmk} One can weaken constraints 5 and 6 on the connection. Namely, it
is enough to assume that the algebra $A$ with product ${\g_0}_{|x}$ 
is unital for every $x\in M$.
\end{rmk}

\subsection{Formal quantization: global $C^{\infty}$-case}

We assume that the reader is familiar with terminology and arguments of
Chapter 7 in \cite{k1}.  Here we shall use a little bit more precise notation
 for the spectrum of the de Rham complex of a manifold, $T[1]\cdot$ instead 
of $\Pi T \cdot$. 

Let $X$ be a $C^{\infty}$-manifold.
 In \cite{k1} we associated with $X$ two differential graded Lie algebras
 $T_{poly}(X)$ and $D_{poly}(X)$, of polyvector fields and of polydifferential
Hochschild cochains on $C^\infty(X)$ respectively, and 
  we proved  that there exists an $L_\infty$ quasi-isomorphism from
 $T_{poly}(X)$ to $D_{poly}(X)$ (see also \cite{cft} for a different approach,
presumably equivalent to one in \cite{k1}).  This quasi-isomorphism depends
on some additional data. Here we propose a more precise result.

 Denote by $Conn_X$ the infinite-dimensional manifold  of all connections
 on the tangent bundle $T_X$. 
\begin{const} There exists a canonical Q-equivariant map
 $$\Phi:T[1]Conn_X\times (T_{poly}(X)[1])_{formal}\ra
(D_{poly}(X)[1])_{formal}$$  such that $\phi(T[1]Conn_X\times \{0\})=\{0\}$
and for any  $\nabla\in Conn_X$ the restriction of  the map $\Phi$  to
$\{\nabla\}\times (T_{poly}(X)[1])_{formal}$  is a quasi-isomorphism. Moreover,
the map $\Phi$ is local (and coordinate-independent) in the obvious sense.
\end{const}

In particular, if $\gamma\in \G(X,\wedge^2 T_X)$ is a Poisson structure,
then the restriction of $\Phi$ to $T[1]Conn_X\times \{\g\}$ gives 
a generalized flat connection (constraints 1-4) on the infinite-dimensional 
affine space  $Conn_X$ with values in the shifted Hochschild complex
$C^\bullet(C^\infty(X),C^\infty(X))[1]$. One can also check  that constraints
5 and 6 from the previous section are satisfied.  Hence the construction of
the previous section applies and we obtain  a canonical algebroid.

The rest of Appendix is devoted to the construction of the map $\Phi$.

{\it Step 1:} With any connection $\nabla$ we associate canonically a
section  of the bundle $X^{aff}\ra X$. Namely, we have a notion of geodesics
for the connection $\nabla$. For very point $x\in X$ we use the formal
completion at $0$ of the exponential map from the tangent space $T_x X$ to 
$X$ in order to get a formal affine structure on $X$ at $x$.
 Applying the functor $T[1]\cdot$ we obtain a $Q$-equivariant map
  $$T[1]Conn_X \ra T[1](Sections(X^{aff}\ra X))$$

{\it Step 2:} We start with a general construction: for any fibration $Y_1\ra
Y_2$ we define a canonical $Q$-equivariant map 
$$T[1](Sections(Y_1\ra Y_2))\ra Sections(T[1]Y_1\ra T[1]Y_2).$$ 
The definition follows easily from the interpretation of  supermanifold 
$T[1]Y$ as the  manifold of maps from $\R^{0|1}$ to $Y$.  
Compatibility with $\Z$-grading and differential $Q$ is a corollary of the 
functoriality (more precisely, it follows from $Diff(\R^{0|1})$-equivariance).
Specifying the above universal construction to the case $Y_1=X^{aff}$ and
$Y_2=X$  we get a map
$$ T[1](Sections(X^{aff}\ra X))\ra Sections(T[1]X^{aff}\ra T[1]X)$$

{\it Step 3:}
Recall that in  subsection 7.3.4 of~\cite{k1} we defined two flat bundles
of differential graded Lie algebras  $\jt$ and $\jd$ on $X$. 
These bundles give two  bundles $(\jt[1])_{formal},\,(\jd[1])_{formal}$ of
formal pointed  graded $Q$-manifolds  over $X$.
 Let us denote by $\pi_1:T[1](X^{aff})\ra X$ the standard projection. 
The result of subsection 7.3.3 of~\cite{k1} can be reformulated as the
construction of $Q$-equivariant mapping  between two bundles of pointed graded
manifolds over $T[1](X^{aff})$ 
$$\pi_1^*((\jt[1])_{formal})\ra\pi_1^*((\jd[1])_{formal})$$
 Let us denote by $\pi_2$ the standard projection of $T[1]X$ to $X$.
 By functoriality  we obtain a new map from the universal map from above:
$$Sections(T[1]X^{aff}\ra T[1]X)\times 
Sections(T[1]X,\pi_2^*(\jt))_{formal}\ra$$  
$$\ra Sections(T[1]X,\pi_2^*(\jd))_{formal}.$$ 
Moreover, this map applied to the submanifold 
$Sections(T[1]X^{aff}\ra T[1]X)\times \{0\}$  gives the point $\{0\}$.

{\it Step 4:} We have a canonical $Q$-equivariant map preserving base points
$$  \alpha_t: T_{poly}(X)[1]_{formal}\ra Sections(T[1]X,\jt)$$
just because $T_{poly}(X)$ is the space of flat sections of $\jt$, 
and an analogous map 
$$ \alpha_d: D_{poly}(X)[1]_{formal}\ra Sections(T[1]X,\jd)$$

These maps are $L_\infty$-quasi-isomorphisms (see 7.3.5. in \cite{k1}).

{\it Quasi-step 5:} Let assume for a moment that we have also a
quasi-isomorphism 
$$ \alpha_d': Sections(T[1]X,\jd)\ra D_{poly}(X)[1]_{formal}.$$ 
Then combining the morphisms listed above we get a promised map.

Unfortunately the situation is more complicated because there is no 
canonical choice of $ \alpha_d'$. One can proceed as follows: 

{\it Step 5:} One can make the following universal construction.
 Let $f:{\bf g}_1\ra {\bf g}_2$  be an $L_\infty$-quasi-isomorphism between
$L_\infty$ algebras ${\bf g}_1$,  ${\bf g}_2$, and we assume 
 that $f_1$ (the first Taylor coefficient of the non-linear map $f$) is a 
monomorphism. Denote by $C_f$ the cokernel of $f_1$. It is a
contractible  complex of vector spaces. Denote by $H_f$
the usual (not graded) manifold of  degree $0$ maps of graded spaces $H:C_f\ra
C_f[1]$ such that $dH+Hd=id_{C_f}$ (the space of homotopies between zero and
the identity map).  Then one can construct a canonical $Q$-equivariant 
map  from $T[1] H_F$ to the supermanifold  of 
$L_\infty$-quasi-isomorphisms  from ${\bf g}_2$ to ${\bf g}_1$. 

We specify this construction to the case  $f= \alpha_d$ and obtain a map 
$$T[1]H_{\alpha_d}\times (D_{poly}(X))_{formal}\ra 
(Sections(T[1]X,\jd))_{formal}.$$

{\it Step 6:} Finally, it is not hard to
construct   a canonical map from the manifold $Sections(X^{aff}\ra X)$ to
$H_{\alpha_d}$. Applying the functor  $T[1]\cdot$ we get a substitute for 
$\alpha_d'$ which IS now parametrized by  $T[1](Sections(X^{aff}\ra X))$, and
combining it with previously constructed maps we obtain the map $\Phi$.


\end{document}